# On the behavior of Bayesian credible intervals for some restricted parameter space problems

Éric Marchand[1] and William E. Strawderman[2]

*Université de Sherbrooke and Rutgers University*

**Abstract:** For estimating a positive normal mean, Zhang and Woodroofe (2003) as well as Roe and Woodroofe (2000) investigate $100(1-\alpha)\%$ HPD credible sets associated with priors obtained as the truncation of noninformative priors onto the restricted parameter space. Namely, they establish the attractive lower bound of $\frac{1-\alpha}{1+\alpha}$ for the frequentist coverage probability of these procedures. In this work, we establish that the lower bound of $\frac{1-\alpha}{1+\alpha}$ is applicable for a substantially more general setting with underlying distributional symmetry, and obtain various other properties. The derivations are unified and are driven by the choice of a right Haar invariant prior. Investigations of non-symmetric models are carried out and similar results are obtained. Namely, (i) we show that the lower bound $\frac{1-\alpha}{1+\alpha}$ still applies for certain types of asymmetry (or skewness), and (ii) we extend results obtained by Zhang and Woodroofe (2002) for estimating the scale parameter of a Fisher distribution; which arises in estimating the ratio of variance components in a one-way balanced random effects ANOVA. Finally, various examples illustrating the wide scope of applications are expanded upon. Examples include estimating parameters in location models and location-scale models, estimating scale parameters in scale models, estimating linear combinations of location parameters such as differences, estimating ratios of scale parameters, and problems with non-independent observations.

## 1. Introduction

For a lower bounded normal mean $\theta$ (say $\theta \geq a$) with unknown standard deviation $\sigma$, and for independent observables $X$ and $W$ with $X \sim \mathrm{N}(\theta, \sigma^2)$, $W \sim \mathrm{Gamma}(\frac{r}{2}, 2\sigma^2)$, Zhang and Woodroofe [9] investigate $100 \times (1-\alpha)\%$ highest posterior density (HPD) credible sets $I_{\pi_0}(X, W)$ associated with the (improper) prior density $\pi_0(\theta, \sigma) = \frac{1}{\sigma} 1_{[a,\infty)}(\theta) 1_{(0,\infty)}(\sigma)$. Using the posterior density $\theta|(X, W)$, which brings into play a truncated Student pdf, they begin by constructing $I_{\pi_0}(X, W)$ as the $100 \times (1-\alpha)\%$ Bayesian interval where the posterior density is the largest.

Then, attractive features of the frequentist coverage of the Bayesian confidence interval $I_{\pi_0}(X, W)$ are established. In particular, they show that

$$(1) \qquad P_{\theta,\sigma}(I_{\pi_0}(X, W) \text{ contains } \theta) \geq \frac{1-\alpha}{1+\alpha},$$

---

[1]Département de mathématiques, Université de Sherbrooke, Sherbrooke Qc, CANADA, J1K 2R1, e-mail: `eric.marchand@usherbrooke.ca`

[2]Department of Statistics, Rutgers University, 501 Hill Center, Busch Campus, Piscataway, NJ 08855, USA, e-mail: `straw@stat.rutgers.edu`

for all $(\theta, \sigma)$ such that $\theta \geq a$ and $\sigma > 0$.

For the case of a known standard deviation, similar developments were given previously by Roe and Woodroofe [6]. Analogously to (1), they obtain that

$$P_\theta(I_{\pi_U}(X) \text{ contains } \theta) \geq \frac{1-\alpha}{1+\alpha}, \tag{2}$$

for all $\theta \geq a$; where $I_{\pi_U}(X)$ is the HPD credible set associated with the prior "uniform" density $\pi_U(\theta) = 1_{[a,\infty)}(\theta)$. Interestingly, for the estimation of the ratio of variance components in a one-way balanced model analysis of variance with random effects, Zhang and Woodroofe obtain [8] results of the same nature.

The objective here is to present extensions of (1) and (2) to other probability models, as well as generalizations to other restricted parameter space scenarios. A notable feature resides in the universal resonance, for symmetric models and for certain types of asymmetric models, of the lower bound $\frac{1-\alpha}{1+\alpha}$. As well, additional frequentist properties of the studied credible intervals are obtained. Although the methods of proof follow for the most part those in the above mentioned papers of Roe, Woodroofe and Zhang, it is particular interesting that the methods of proof are unified. Moreover, we actually offer a useful simplification.

Inference problems for constrained parameter spaces has, for many years, held the interest of statisticians. Correspondingly, as reviewed by Marchand and Strawderman [5] or van Eeden [7], it has been a fairly active field. Recently though, there has been a renewed interest from the particle physicist community with high energy experiments leading to constrained parameter models (see for instance [2, 3, 4]), and more specifically to the problem of setting confidence bounds in the presence of constrained parameters. Actually, a vigorous and substantial debate has arisen, focussing indeed on the choice of method, with an underlying Bayesian-frequentist comparison of the respective advantages and disadvantages (e.g., [4]).

As an example for the normal model above with known variance, it has been observed that the so-called "unified method" put forth by Feldman and Cousins [3]; which is a frequentist based method arrived at by the inversion of a LRT and which leads to exact frequentist coverage; produces "quite short" intervals for small values of $X$, in comparison at least to $I_{\pi_U}(X)$. Such observations are not surprising since the methods differ in how they take into account the lower-bound constraint. As argued by Zhang and Woodroofe [9] and Roe and Woodroofe [6], the HPD credible intervals $I_{\pi_0}(X, W)$ and $I_{\pi_U}(X)$ are quite sensible ways to deal with the lower bound constraint. If such is the case, then good frequentist coverage properties of these Bayesian confidence intervals would render them more attractive, even from a frequentist point of view. There lies as well an intrinsic interest in these procedures given that the untruncated versions of the priors $\pi_0$ and $\pi_U$ lead to the usual, and introductory textbook, $t$ and $z$ two-sided $100 \times (1-\alpha)\%$ intervals; which of course have exact coverage.

The paper is organized as follows. Symmetric models are treated in Section 2, while asymmetric models are reserved for Section 4. The main finding of Section 2 relates to the choice of the truncation of the Haar right invariant prior for a large class of problems, with underlying unimodality and symmetry, which leads to the lower bound $\frac{1-\alpha}{1+\alpha}$ for the frequentist coverage probability of the associated HPD credible set. Various other corollaries are available. For instance, an exact coverage probability of $\frac{1}{1+\alpha}$ for boundary parameter values is established, and the above $\frac{1-\alpha}{1+\alpha}$ lower bound is shown to hold for a Bayesian confidence interval which is not HPD (see Remark 1, part c). Various examples, illustrating the wide scope of applications, are expanded upon in Section 3. The developments for asymmetric models is more delicate requiring a categorization of different types of asymmetry. In cases where the underlying models' density is not monotone, the lower bounds obtained in general are less explicit, but there is evidence that these lower bounds can be quite large. Moreover, the last result (Corollary 2) actually recovers the lower bound $\frac{1-\alpha}{1+\alpha}$ for certain types of underlying skewness, as a generalization of the symmetric case.

## 2. Symmetric models

We first consider models with an observable scalar or vector $X$ having densities $f(x;\theta)$; $\theta \in A \subset \Re^p$; for which there exists a lower bound constraint of the form $\tau(\theta) \geq 0$; $\tau(\theta) : \Re^p \to \Re$. Moreover, we work with a structure, which is present in previous work described above, and where there exists a linear pivotal quantity of the form $\frac{a_1(X)-\tau(\theta)}{a_2(X)}$ with underlying absolutely continuous, symmetric and (strictly) unimodal density. An immediate example consists of symmetric and unimodal location densities $f(x;\theta) = f_0(x-\theta)$, with $\tau(\theta) = \theta \geq 0$ and the pivot $X - \theta$. Further examples are presented in Section 3.

We study HPD credible intervals $I_{\pi_0}(X)$, based on $(a_1(X), a_2(X))$, associated with a prior $\pi_0$ obtained as the truncation onto the parameter space $\{\theta : \tau(\theta) \geq 0\}$ of a Haar right invariant density $\pi(\theta)$. To describe the construction of this interval as well as several others that follow, it is useful to define the following quantities.

**Definition 1.** For a given continuous cumulative distribution function $F$, $\alpha \in (0,1)$, $y \in \Re$, we define:

$$d_1(y) = F^{-1}(1 - \alpha F(y)); d_2(y) = F^{-1}(\frac{1}{2} + \frac{1-\alpha}{2} F(y)); \text{ and}$$
$$d(y) = \max(d_1(y), d_2(y)).$$

In situations where we wish to emphasize the dependence of the above functions on the pair $(F, \alpha)$, we will write instead $d_{1_{F,\alpha}}$, $d_{2_{F,\alpha}}$, and $d_{F,\alpha}$. To a large extent, the frequentist properties which we establish below depend on the following property of $d_{F,\alpha}$; which is easily established.

**Lemma 1.** For all $(F, \alpha)$, we have $d(y) \geq d(d_0) = d_0$, with $d_0 = F^{-1}(\frac{1}{1+\alpha})$. As well, $d(y) = d_1(y)$ if and only if $y \leq d_0$.

*Proof.* A direct evaluation tells us that $d_0 = F^{-1}(\frac{1}{1+\alpha})$ is a fixed point of $d_1, d_2$, and hence of $d$. The result follows as $d_1(y)$ decreases in $y$, and $d_2(y)$ increases in $y$. □

The following theorem is our first key result. Paired with Corollary 1, it will lead to various applications which are generalizations of (1) and (2).

**Theorem 1.** *For a model $X|\theta \sim f(x;\theta)$, and a parametric function $\tau(\theta) : \Re^p \to \Re$ such that $\tau(\theta) \geq 0$ (constraint); suppose there exists a linear pivot $T(X, \theta) = \frac{a_1(X)-\tau(\theta)}{a_2(X)}$; with $a_2(\cdot) > 0$; such that the distribution of $T(X, \theta)|\theta$ is given by cdf $G$, with pdf $G'$ which is symmetric and unimodal (without loss of generality, about 0). Suppose further that there exists a prior $\pi(\theta)$ supported on the natural parameter space such that:*

(3) $$T(X,\theta)|x =^d T(X,\theta)|\theta;$$

(i.e., the frequentist distribution of $T(X,\theta)$ for a given $\theta$; which is independent of $\theta$ and given by cdf $G$; matches the posterior distribution of $T(X,\theta)$ for any given value $x$ of $X$). Then, for the prior $\pi_0(\theta) = \pi(\theta)1_{[0,\infty)}(\tau(\theta))$, we have:

(a) $I_{\pi_0}(X) = [l(X), u(X)]$, with $l(X) = \max\{0, a_1(X) - d_{G,\alpha}(\frac{a_1(X)}{a_2(X)})a_2(X)\}$ and $u(X) = a_1(X) + d_{G,\alpha}(\frac{a_1(X)}{a_2(X)})a_2(X)$;
(b) $P_\theta(I_{\pi_0}(X) \ni \tau(\theta)) > \frac{1-\alpha}{1+\alpha}$, for all $\theta$ such that $\tau(\theta) \geq 0$;
(c) $P_\theta(I_{\pi_0}(X) \ni 0) = \frac{1}{1+\alpha}$, for all $\theta$ such that $\tau(\theta) = 0$;
(d) $\lim_{\tau(\theta) \to \infty} P_\theta(I_{\pi_0}(X) \ni \tau(\theta)) = 1 - \alpha$.[1]

*Proof.* **(a)** Denote $h_x$, $H_x$, and $H_x^{-1}$ as the pdf, cdf, and inverse cdf of the posterior distribution of $\tau(\theta)$ under $\pi_0$. Since $T(X,\theta)$ is a pivot, implying that its distribution is, for any given $\theta$, free of $\theta$, we infer from (3) that, for $\theta \sim \pi$, $P_\pi(T(X,\theta) \leq y|x) = G(y)$, or equivalently $P_\pi(\tau(\theta) \geq y|x) = G(\frac{a_1(x)-y}{a_2(x)})$. By definition of $\pi_0$, this gives us for $y \geq 0$

$$H_x(y) = P_{\pi_0}(\tau(\theta) \leq y|x) = 1 - \frac{P_\pi(\tau(\theta) \geq y|x)}{P_\pi(\tau(\theta) \geq 0|x)} = 1 - \frac{G(\frac{a_1(x)-y}{a_2(x)})}{G(\frac{a_1(x)}{a_2(x)})};$$

and

$$H_x^{-1}(\Delta) = a_1(x) - a_2(x)G^{-1}((1-\Delta)G(\frac{a_1(x)}{a_2(x)})).$$

Now, observe that the posterior density is unimodal, with a maximum at $\max(0, a_1(x))$. From this, since our HPD credible interval may be represented as $\{\tau(\theta) : h_x(\tau(\theta)) \geq c\}$ for some constant $c$ (e.g., [1], page 140), we infer that either:

(i) $l(x) = 0$ and $u(x) = H_x^{-1}(1-\alpha) = a_1(x) - a_2(x)G^{-1}(\alpha G(\frac{a_1(x)}{a_2(x)}))$, or
(ii) $l(x) = a_1(x) - b(x)$ and $u(x) = a_1(x) + b(x)$; for some $b(x)$ such that $a_1(x) - b(x) > 0$.

From the symmetry of $G'$, we have in (i): $u(x) = a_1(x) + a_2(x)G^{-1} \times (1 - \alpha G(\frac{a_1(x)}{a_2(x)})) = a_1(x) + a_2(x) d_{1_{G,\alpha}}(\frac{a_1(x)}{a_2(x)})$. For (ii), we obtain also with the symmetry of $G'$ that:

$$P_{\pi_0}(a_1(x) - b(x) \leq \tau(\theta) \leq a_1(x) + b(x)|x) = 1 - \alpha$$
$$\Leftrightarrow H_x(a_1(x) + b(x)) - H_x(a_1(x) - b(x)) = 1 - \alpha$$
$$\Leftrightarrow G(\frac{b(x)}{a_2(x)}) - G(-\frac{b(x)}{a_2(x)}) = (1-\alpha)G(\frac{a_1(x)}{a_2(x)})$$
$$\Leftrightarrow G(\frac{b(x)}{a_2(x)}) = \frac{1}{2}[1 + (1-\alpha)G(\frac{a_1(x)}{a_2(x)})] \text{ (by symmetry)}$$
$$\Leftrightarrow b(x) = a_2(x)G^{-1}(\frac{1}{2} + \frac{(1-\alpha)}{2}G(\frac{a_1(x)}{a_2(x)})) = a_2(x) d_{2_{G,\alpha}}(\frac{a_1(x)}{a_2(x)}).$$

Moreover, situation (ii) occurs iff

(4) $$a_1(x) > a_2(x) d_{2_{G,\alpha}}(\frac{a_1(x)}{a_2(x)}) \Longleftrightarrow G(\frac{a_1(x)}{a_2(x)}) > \frac{1}{2} + \frac{1-\alpha}{2}G(\frac{a_1(x)}{a_2(x)})$$
$$\Longleftrightarrow \frac{a_1(x)}{a_2(x)} > G^{-1}(\frac{1}{1+\alpha}) = d_0.$$

---

[1] More precisely, we are referring of course to a sequence of $\theta_i's; i = 1, 2, \ldots$ such the corresponding $\tau(\theta_i)$'s have a limiting value of $+\infty$.

Finally, the result follows by combining (i) and (ii) and using Lemma 1.

**(b)** First, observe that the interval $a_1(X) \pm d_{G,\alpha}(\frac{a_1(X)}{a_2(X)}) \, a_2(X)$ has the same coverage probability as its subset $I_{\pi_0}(X)$ for nonnegative values of $\tau(\theta)$, since the difference of these two sets can only help in covering negative values of $\tau(\theta)$. Now, along the lower bound $d_0$ of Lemma 1 and the symmetry of $G'$, we have for $\theta$'s such that $\tau(\theta) \geq 0$:

$$P_\theta(I_{\pi_0}(X) \ni \tau(\theta)) = P_\theta(a_1(X) - a_2(X)d_{G,\alpha}(\frac{a_1(X)}{a_2(X)})$$
$$\leq \tau(\theta) \leq a_1(X) + a_2(X)d_{G,\alpha}(\frac{a_1(X)}{a_2(X)}))$$
$$= P_\theta(|T(X,\theta)| \leq d_{G,\alpha}(\frac{a_1(X)}{a_2(X)}))$$
$$> P_\theta(|T(X,\theta)| \leq d_0)$$
$$= 2G(d_0) - 1 = 2G(G^{-1}(\frac{1}{1+\alpha})) - 1 = \frac{1-\alpha}{1+\alpha}.$$

**(c)** Since coverage at $\tau(\theta) = 0$ occurs if and only $l(X) = 0$, we have by (4) for $\theta$'s such that $\tau(\theta) = 0$:

$$P_\theta(I_{\pi_0}(X) \ni 0) = P_\theta(l(X) = 0) = P_\theta(\frac{a_1(X) - 0}{a_2(X)} \leq d_0) = G(d_0) = \frac{1}{1+\alpha}.$$

**(d)** Since $\frac{a_1(X) - \tau(\theta)}{a_2(X)}$ is a pivot, implying that "$\frac{a_1(X)}{a_2(X)} \to \infty$" and $G(\frac{a_1(X)}{a_2(X)})$ converges to 1 in probability as $\tau(\theta) \to \infty$, it follows that $d_{G,\alpha}(\frac{a_1(X)}{a_2(X)})$ (equal to $G^{-1}(\frac{1}{2} + \frac{1-\alpha}{2}G(\frac{a_1(X)}{a_2(X)}))$ for large $\frac{a_1(X)}{a_2(X)}$) converges in probability, as $\tau(\theta) \to \infty$, to $G^{-1}(1 - \frac{\alpha}{2})$. In view of the above, and as in part (b), we have

$$\lim_{\tau(\theta)\to\infty} P_\theta(I_{\pi_0}(X) \ni \tau(\theta)) = \lim_{\tau(\theta)\to\infty} P_\theta(|T(X,\theta)| \leq d_{G,\alpha}(\frac{a_1(X)}{a_2(X)}))$$
$$= P_\theta(|T(X,\theta)| \leq G^{-1}(1 - \frac{\alpha}{2}))$$
$$= 2G(G^{-1}(1 - \frac{\alpha}{2})) - 1 = 1 - \alpha. \qquad \square$$

Observe how critical (3) is, namely in the last line of the proof of part (b) where the identity $G(G^{-1})$ arises. In fact, the "$G^{-1}$" comes from the construction of $I_{\pi_0}(X)$ (hence the lhs of (3)), while the "$G$" comes from the frequentist coverage assessment of $I_{\pi_0}(X)$ (hence the rhs of (3)). Condition (3) may appear stringent but, as shown below, it is attainable for a large class of problems if the prior $\pi(\theta)$ is Haar right invariant (informally, a prior leaving the measure of sets constant under certain transformations). For instance, consider a simple location model $X \sim f_0(x-\theta)$ with known $f_0$. Set $Z = X - \theta$ and consider the flat prior $\pi(\theta) = 1$. It is easy to verify that for any pair $(x,\theta)$, the distributions $Z|\theta$ and $Z|x$ match with density $f_0(\cdot)$, which tells us that condition (3) holds here with the choice of the flat (right Haar invariant also) prior. (It is important to note that the assumptions of symmetry and unimodality are additional and specific to Theorem 1, and are not required for the above illustration of (3). This is exploited namely in Section 4 (also see Remark 1, part c) where we make use of condition (3)).

The various applications (see Section 3) which will follow from Theorem 1 are essentially all cases where the prior $\pi(\theta)$ is Haar right invariant (denoted $\pi^r(\theta)$) and the pivot satisfies the invariance requirement $T(x,\theta) = T(gx, \bar{g}\theta)$, for all $x \in \mathcal{X}$, $\theta \in \Theta$, $g \in G$, $\bar{g} \in \bar{G}$, with $\mathcal{X}$, $\Theta$, $G$, and $\bar{G}$ being isomorphic ("equivalent"). We now pursue by showing how this invariance requirement and conditions lead to (3), hence permitting the application of Theorem 1 for a given problem. We make use of the following result (and notation) given in [1].

**Lemma 2** (Result 3, p. 410 [1]). *Consider an invariant decision problem for which $\mathcal{X}$, $\Theta$, $G$, and $\bar{G}$ are all isomorphic. Then, for an invariant decision rule $\delta(x) = \tilde{x}(a)$,*

$$(5) \qquad E^{\pi^r(\theta|x)}\{L(\theta, \tilde{x}(a))\} = R(\theta, \delta(X)),$$

*where $\pi^r(\theta|x)$ is the posterior distribution with respect to the right invariant (generalized) prior density $\pi^r(\theta)$.*

**Corollary 1.** *Suppose $\mathcal{X}$, $\Theta$, $G$, and $\bar{G}$ are all isomorphic, and that $T(X, \theta)$ is a function for which $T(x, \theta) = T(gx, \bar{g}\theta)$, for all $x \in \mathcal{X}$, $\theta \in \Theta$, $g \in G$, $\bar{g} \in \bar{G}$. Then condition (3) holds, that is $P_\theta[T(X,\theta) \in A] = P^{\pi^r(\theta|x)}[T(X,\theta) \in A]$ for each measurable set $A$ (where the lhs gives the frequentist distribution of $T(X,\theta)$ for given $\theta$, and the rhs gives the posterior distribution of $T(X,\theta)$ conditional on $X$ for the right invariant Haar measure).*

*Proof.* It suffices to establish, for each measurable set $A$ (in the range of $T(X,\theta)$), the identity:

$$(6) \qquad P_\theta(T(X,\theta) \in A) = P^{\pi^r(\theta|x)}(T(X,\theta) \in A).$$

To do so, we apply Lemma 2 for loss $L_A(\theta, d) = 1_A(T(d, \theta))$, and for $\delta(X) = X$. With $G = G^*$, we indeed have that $\delta(X)$ is an equivariant decision rule since $\delta(gX) = gX = g^*(X) = g^*(\delta(X))$. We also have by assumption on $T$:

$$L_A(\bar{g}\theta, g^*d) = L_A(\bar{g}\theta, gd) = 1_A(T(gd, \bar{g}\theta)) = 1_A(T(d, \theta)) = L_A(\theta, d),$$

which tells us that we have an invariant decision problem. Finally, applying Lemma 2 yields (6) and establishes the Corollary. □

Now, prior to presenting various illustrations and applications of Theorem 1 (and Corollary 1) in Section 3, we conclude this section by expanding on some interesting aspects and implications of the results above.

**Remark 1.**  (a) Exact values or very good approximations of the frequentist coverage probability of $I_{\pi_0}$, which seem difficult to establish, are not provided explicitly by the results above. The exceptions are at the boundary where the probability of coverage $\frac{1}{1+\alpha}$ exceeds the nominal coverage probability $1-\alpha$, and when $\tau(\theta) \to \infty$ where the coverage probability tends to $1-\alpha$. Numerical evaluations are provided, for the normal models described in the introduction, by Roe and Woodroofe [6], and Zhang and Woodroofe [9]. Moreover, as pointed out in these manuscripts for a normal model $G$, and as suggested by the derivation above, the lower bound $\frac{1-\alpha}{1+\alpha}$ is, for a specific $G$, somewhat conservative. But it has the advantage of being simple and derived in a unified fashion, applicable for a vast array of situations, and for quite general symmetric and unimodal densities $G'$.
 (b) In addition, the above development can be adapted to deal with the following robustness issue. Indeed, suppose that the actual model is governed by symmetric pdf's $f_1(x; \theta)$, with corresponding cdf's $G_1$, in contrast to the bounds which are set using $G$. Then, following the proof of Theorem 1(b), above, we have

$$P_\theta(I_{\pi_0}(X) \ni \tau(\theta)) > P_0(|T(X,0)| \leq d_0) = 2G_1(d_0) - 1;$$

which provides lower bounds or envelopes depending on $G_1$. Moreover, the quantity $\frac{1-\alpha}{1+\alpha}$ remains a lower bound on the probability of coverage for a given $G_1$ as long as, simply,

(7) $$G_1(d_0) \geq G(d_0); \text{ (with } d_0 = G^{-1}(\frac{1}{1+\alpha})).$$

Here, various properties of families of distributions can be elucidated to give realizations of (7). For instance, (7) holds as long as $\frac{G'_1(y)}{G'(y)}$ is nonincreasing in $y; y > 0$; or as a specific case if $G'_1(y) = \frac{1}{\sigma}G'(\frac{y}{\sigma})$, $\sigma < 1$, in other words $f_1$ and $f_0$ belong to the same scale family having increasing monotone likelihood ratio in $|y|$.

(c) Interestingly, in the case of continuous but non-unimodal $G'$, the above development remains valid with the difference that the interval $I_{\pi_0}(X)$ is not HPD, in other words $I_{\pi_0}(X)$ is a credible interval with the same frequentist properties as those given in Theorem 1, but it is not (necessarily) optimal in the sense of being the credible region with the shortest length.

## 3. Examples

We enumerate a list of situations for which Theorem 1 applies. The list is also illustrative in the sense that we also specify components, such as the pivot $T(X,\theta)$ and the prior $\pi_0$ of Theorem 1. In all cases below with unimodal and symmetric density $G'$, the lower bound $\frac{1-\alpha}{1+\alpha}$ applies for the coverage probability of the confidence interval $I_{\pi_0}(X)$. In cases where the density $G'$ is unimodal but not symmetric, the results of Section 4 will also apply to each one of the following situations as well.

(a) (location) $X \sim f_0(x-\theta)$; $f_0$ unimodal and symmetric; $\tau(\theta) = \theta \geq 0$; $T(X,\theta) = X - \theta$; $\pi_0(\theta) = 1_{[0,\infty)}(\theta)$. For example, this applies for a $N(\theta,\sigma)$ model known $\sigma$ and $\theta \geq 0$; but also to many other common univariate symmetric models such as Logistic, Laplace, Cauchy and Student, etc.

(b) (location-scale) $(X_1, X_2) \sim f_0(\frac{x_1-\theta_1}{\theta_2}, \frac{x_2}{\theta_2})$; $\tau(\theta) = \theta_1 \geq 0$; $T(X,\theta) = \frac{X_1-\theta_1}{X_2}$; $\pi_0(\theta) = \frac{1}{\theta_2}1_{(0,\infty)}(\theta_2)1_{[0,\infty)}(\theta_1)$. Observe that $T(X,\theta)$ is indeed a pivot here as it can expressed as the ratio of the elements of the pair $(\frac{X_1-\theta_1}{\theta_2}, \frac{X_2}{\theta_2})$, whose distribution is free of $(\theta_1,\theta_2)$. An important case here arises with the model $Y_1,\ldots,Y_n$ i.i.d. $N(\theta_1,(\theta_2)^2)$, and for which the sufficient statistic $(X_1, X_2) = (\bar{Y}, \frac{S_y}{\sqrt{n}})$ admits a location-scale model as above with the distribution of $T(X,\theta)|\theta$ being Student with $n-1$ degrees of freedom.

(c) (multivariate location) $X = (X_1,\ldots,X_p) \sim f_0(x_1-\theta_1,\ldots,x_p-\theta_p)$; $\tau(\theta) = \sum_{i=1}^p a_i\theta_i$; $T(X,\theta) = (\sum_{i=1}^p a_iX_i) - \tau(\theta)$; $\pi_0(\theta) = 1_{[0,\infty)}(\tau(\theta))$. For example, take $X \sim N_p(\theta,\Sigma)$; $\Sigma$ known; in which case $T(X,\theta)|\theta \sim N(0, a'\Sigma a)$, with $a' = (a_1,\ldots,a_p)$. An important case here (and in (d) as well) concerns the estimation of the difference of two means $\theta_1 - \theta_2$, with the information that $\theta_1 \geq \theta_2$.

(d) (multivariate location-scale with homogeneous scale)
$X = (X_1,\ldots,X_p,X_{p+1}) \sim f_0(\frac{x_1-\theta_1}{\theta_{p+1}},\ldots,\frac{x_p-\theta_p}{\theta_{p+1}},\frac{x_{p+1}}{\theta_{p+1}})$; $\tau(\theta) = \sum_{i=1}^p a_i\theta_i$;

$T(X, \theta) = \frac{(\sum_{i=1}^{p} a_i X_i) - \tau(\theta)}{X_{p+1}}$; $\pi_0(\theta) = \frac{1}{\theta_{p+1}} 1_{(0,\infty)}(\theta_{p+1}) 1_{[0,\infty)}(\tau(\theta))$. For example, consider $(X_1, \ldots, X_p)'$ and $X_{p+1}$ independent with $(X_1, \ldots, X_p) \sim N_p((\theta_1, \ldots, \theta_p), \theta_{p+1}^2 I_p)$ and $X_{p+1}^2 \sim \text{Gamma}(r/2, 2\theta_{p+1}^2)$, in which case the distribution of $T(X, \theta)|\theta$ is distributed as $\{(\sum_{i=1}^n a_i^2)\sqrt{\frac{2}{r}}\} T_r$, with $T_r$ distributed Student with $r$ degrees of freedom.

(e) (scale with support being a subset of $\Re^+$ or $\Re^+$) $X \sim \frac{1}{\theta} f_1(\frac{x}{\theta})$; $\tau(\theta) = \log(\theta) - \log(a) \geq 0$; $T(X, \theta) = \log(\frac{X}{\theta})$ (i.e., $a_1(X) = \log(X) - \log(a), a_2(X) = 1$); $\pi_0(\theta) = \frac{1}{\theta} 1_{[0,\infty)}(\tau(\theta))$. The constraint on $\tau(\theta)$ corresponds to a lower bound constraint on $\theta$, and confidence intervals for $\tau(\theta)$ provide confidence intervals for $\theta$, with corresponding frequentist coverage probabilities. As a specific example, consider a lognormal model with scale parameter $\theta$; $\theta \geq a(> 0)$; where $\frac{X}{\theta} \sim e^{\delta Z}$, $Z \sim N(0,1)$, and $\delta$ being a known and positive shape parameter. Here $f_1(y) = (\sqrt{2\pi}\delta y)^{-1} e^{-\frac{(\log y)^2}{2\delta^2}} 1_{(0,\infty)}(y)$, and the distribution of $T(X, \theta)|\theta$ is normal with mean 0 and standard deviation $\delta$. Additional examples arise from scale models such that $\frac{X}{\theta}$ and $(\frac{X}{\theta})^{-1}$ are equidistributed which implies symmetry for the distribution of the pivot $T(X, \theta) = \log(\frac{X}{\theta})$. Further specific examples where $\frac{X}{\theta}$ and $(\frac{X}{\theta})^{-1}$ are equidistributed include the half-Cauchy with $f_1(y) = \frac{2}{\pi}(1+y^2)^{-1} 1_{(0,\infty)}(y)$, and Fisher distributions with matching degrees of freedom in both numerator and denominator. On the other hand, if $X \sim \text{Gamma}(\alpha, \theta)$ for instance, then the distribution of $\log(\frac{X}{\theta})$ is not symmetric (for any $\alpha$), but the results of Section 4 apply nevertheless (see Example 2).

(f) (multivariate scale) $(X_1, \ldots, X_p) \sim (\Pi_{i=1}^p \frac{1}{\theta_i}) f_1(\frac{x_1}{\theta_1}, \ldots, \frac{x_p}{\theta_p})$; $\tau(\theta) = \sum_{i=1}^p a_i \log(\theta_i)$. For instance in correspondence to the problem of estimating the ratio of two scale parameters under the lower bound constraint $\frac{\theta_2}{\theta_1} \geq a$, $\tau(\theta) = \log(\theta_2) - \log(\theta_1) - \log(a)$; $T(X, \theta) = \log(\frac{X_2}{\theta_2}) - \log(\frac{X_1}{\theta_1}) - \log(a)$ (i.e, $a_1(X) = \log(\frac{X_2}{X_1}) - \log(a), a_2(X) = 1$); $\pi(\theta) = \frac{1}{\theta_1 \theta_2} 1_{[a,\infty)}(\frac{\theta_2}{\theta_1})$. Specific examples here arise whenever $X_1$ and $X_2$ are independent with the distributions of $\log(\frac{X_i}{\theta_i})$; $i = 1, 2$; being symmetric (see part (e) above). Hence, Theorem 1 can be applied for instance to estimating a lower-bounded ratio of two lognormal scale parameters.

We note that none of the above situations requires independence between the vector components (and see Example (g)). Theorem 1 applied to Example (a) and (b) extends the results of Roe and Woodroofe [6], Zhang and Woodroofe [8, 9] obtained for the normal case. The asymmetric case studied by Zhang and Woodroofe (2002) which deals with a Fisher distribution is contained in part (e) (here Theorem 2 and perhaps Corollary 2 apply). Numerical displays of $I_{\pi_0}(X)$ and of its coverage probability, in comparison namely to other confidence interval procedures, are given in the above papers, as well in [4].

The developments above are neither limited to samples of size 1 of $X$, nor to cases where $X$ is a sufficient statistic. Further applications are available by conditioning on the maximal invariant $V$. For instance, the results are applicable for location parameter families with densities $f(x_1 - \theta, \ldots, x_n - \theta)$, provided the conditional distribution of $\bar{X}_n = \sum_{i=1}^n X_i$ given the maximal invariant $v = (x_1 - \bar{x}_n, \ldots, x_{n-1} - \bar{x}_n)$ satisfies the conditions required for $G$ (a.e. $v$). We conclude this section with an illustration with spherically symmetric models, and specifically to a multivariate student model.

(g) (sample of size $n$ with underlying spherically symmetric distribution) Suppose the distribution of $X = (X_1, \ldots, X_n)$ is spherically symmetric about $(\theta, \ldots, \theta)$ with density $f(x; \theta) = h(\sum_{i=1}^n (x_i - \theta)^2)$, or equivalently,

$$f(x;\theta) = h(n(\bar{x}_n - \theta)^2 + \sum_{i=1}^n (x_i - \bar{x}_n)^2). \tag{8}$$

Considering now the pivot $Z = T(X, \theta) = \bar{X}_n - \theta$ and the maximal invariant $V = (X_1 - \bar{X}_n, \ldots, X_{n-1} - \bar{X}_n)$, Theorem 1 applies for the procedure $I_{\pi_0}(X, V)$ which is constructed as in part (a) of Theorem 1 but with the cdfs $G_v$ associated with the conditional distributions $Z|V = v$, or equivalently by virtue of (8) with the conditional pdfs

$$G'_{Z|v}(z) \propto h(nz^2 + B(v)); \tag{9}$$

with $B(v) = v'(I_{n-1} + 11')v$, $1' = (1, \ldots, 1)$; as $B(v) = \sum_{i=1}^n (x_i - \bar{x}_n)^2$. The key points being that the conditional distributions $Z|V = v$ are free of $\theta$, and that the bounds on conditional coverage associated with $G_v$ are free of $v$. As a specific example, consider a multivariate Student model for $X = (X_1, \ldots, X_n)$ with $d$ degrees of freedom, location parameter $(\theta, \ldots, \theta)$, scale parameter $\sigma$, such that $h(y) \propto (1 + \frac{y}{d\sigma^2})^{-(\frac{d+n}{2})}$ in (8). An evaluation of (9) tells us that $G'_{Z|v}(z) \propto (1 + \frac{nz^2 + B(v)}{d\sigma^2})^{-(\frac{d+n}{2})} \propto (1 + \frac{z^2}{\nu \sigma'^2})^{-(\frac{\nu+1}{2})}$, with $\nu = d + n - 1$ and $\sigma'^2 = \frac{\sigma^2 d + B(v)}{n(d+n-1)}$. In other words, the conditional cdfs $G_{Z|v}$, which are used to construct $I_{\pi_0}(X, v)$, are those of a univariate Student distribution with degrees of freedom $d + n - 1$ and scale parameter $\sigma' = \sqrt{\frac{\sigma^2 d + B(v)}{n(d+n-1)}}$.

## 4. Asymmetric models

Here, we investigate and extend the results of Section 2 to unimodal, but not necessarily symmetric densities. However, as illustrated with the next example, unified lower bounds on the frequentist coverage probability, such as those given in Theorem 1, are not possible and conditions on the type of asymmetry are required.

**Example 1.** Consider an exponential location model with density $e^{-(x-\theta)} 1_{(0,\infty)}(x - \theta)$; and $\theta \geq 0$. For the uniform prior $\pi_0(\theta) = I_{[0,\infty)}(\theta)$, the $(1-\alpha) \times 100\%$ HPD credible interval is given by $I_{\pi_0}(X) = [l(X), u(X)]$, with $l(x) = \log(1 - \alpha + \alpha e^x)$ and $u(x) = x$. Observe that the interval never covers the value $\theta = 0$, so that the coverage probability $P_0(I_{\pi_0}(X) \ni 0)$ is equal to 0. Hence, a very different situation arises in comparison to the case of symmetric $G's$. Moreover, it is easy to establish that

$$\begin{aligned} P_\theta(I_{\pi_0}(X) \ni \theta) &= P_\theta(\theta \geq \log(1 - \alpha + \alpha e^X)) \\ &= P_\theta(X \leq \log(1 + \frac{e^\theta - 1}{\alpha})) \\ &= 1 - e^{-(\log(1 + \frac{e^\theta - 1}{\alpha}) - \theta)} \\ &= (1 - \alpha) \frac{e^\theta - 1}{\alpha + e^\theta - 1}. \end{aligned}$$

Hence, the coverage probability can be quite small and never exceeds the nominal coverage level $1 - \alpha$.[2] Finally, as one may anticipate, the same characteristics will arise for more general models with a property of monotone decreasing densities (see Theorem 2, part b).

As in Theorem 1 and Corollary 1, the results below apply to models $X|\theta \sim f(x; \theta)$ and for estimating $\tau(\theta)$ under the constraint $\tau(\theta) \geq 0$.

**Assumption 1.** We assume again the existence of a linear pivot $T(X, \theta) = \frac{a_1(X) - \tau(\theta)}{a_2(X)}$ such that $-T(X, \theta)$ has cdf $G$, with (strict) unimodal $G'$. Moreover, we assume without loss of generality that the density $G'$ has a mode at 0.

The confidence interval procedures studied are HPD credible based on $(a_1(X), a_2(X))$, and associated with the truncation $\pi_0$ of the Haar right-invariant $\pi^r$ onto the constrained parameter space; i.e, $\pi_0(\theta) = \pi^r(\theta) I_{[0,\infty)}(\tau(\theta))$. We pursue with the introduction of various quantities and related properties which will help in describing the $(1 - \alpha) \times 100\%$ HPD credible interval $I_{\pi_0}(X)$, as well as some of its frequentist properties. In particular, as illustrated above in the contrasting results of Example 1 and of Theorem 1, and since the frequentist properties which we can hope to establish depend on the type of asymmetry present, we breakdown, in Definition 4 and Corollary 2, these asymmetries into different relevant types. This is achieved in part with the introduction of the function $U_{G,\alpha}$ in Definition 3 below; which will also relate to familiar qualitative features such as skewness to the right (see Corollary 2).

**Definition 2.** For cdf $G$ with unimodal at 0 density $G'$, and $\Delta \in (0, 1)$, define $\gamma_1(\Delta)$ and $\gamma_2(\Delta)$ as values that minimize the length $|\gamma_1 + \gamma_2|$ among all intervals $[-\gamma_1, \gamma_2]$ such that $G(\gamma_2) - G(-\gamma_1) = \Delta$.

Observe that the above defined $\gamma_1(\Delta)$ and $\gamma_2(\Delta)$ are indeed uniquely determined, and nonnegative given the unimodality. Furthermore, note that if $G(0) \in (0, 1)$, then we also have $G'(-\gamma_1(\Delta)) = G'(\gamma_2(\Delta))$.

**Definition 3.** Let $1 - \alpha \in (0, 1)$ and $G$ be a cdf with unimodal density $G'$ with a mode at 0. Let

$$(10) \qquad U_{G,\alpha}(y) = -y + \gamma_1((1 - \alpha)(1 - G(-y)));$$

be defined for values $y$ such that $-y$ belongs to the support of $G'$, (i.e., $y \in (-G^{-1}(1), -G^{-1}(0))$).

**Definition 4.** Let $1 - \alpha \in (0, 1)$. Let $\mathcal{C}_1$, $\mathcal{C}_2$, and $\mathcal{C}_3$ be classes of cdfs $G$ with unimodal at 0 densities $G'$ such that

$$\mathcal{C}_1 = \{G : \text{there exists an interior point } y_0 \in (-G^{-1}(1), -G^{-1}(0)) \\ \text{such that } U_{G,\alpha}(y_0) = 0\}$$

$\mathcal{C}_2 = \{G : U_{G,\alpha}(y) \geq 0 \text{ for all -y on the support of G'}\}$

$\mathcal{C}_3 = \{G : U_{G,\alpha}(y) \leq 0 \text{ for all -y on the support of G'}\}$

**Lemma 3.** *In the context of Definition 4, the classes $\mathcal{C}_1$, $\mathcal{C}_2$, and $\mathcal{C}_3$ can alternatively be described as $\mathcal{C}_1 = \{G : G(0) \in (0, 1)\}, \mathcal{C}_2 = \{G : G(0) = 0\}$, and $\mathcal{C}_3 = \{G : G(0) = 1\}$.*

---

[2] On the other hand, the coverage rises fast as $\theta$ increases and attains, for instance, Theorem 1's lower bound $\frac{1-\alpha}{1+\alpha}$ as soon as $\theta = \log 2$.

**Note.** In other words, the class $\mathcal{C}_2$ consists of decreasing densities $G'$; the class $\mathcal{C}_3$ consists of increasing densities $G'$, and $\mathcal{C}_1$ consists of densities $G'$ which increase on $\Re^-$ and decrease on $\Re^+$.

*Proof.* First observe that

$$U_{G,\alpha}(y)|_{y=-G^{-1}(0)} = G^{-1}(0) + (\gamma_1(1-\alpha)) \leq 0,$$

with equality iff $G^{-1}(0) = \gamma_1(1-\alpha) = 0$, i.e., $G(0) = 0$. Similarly,

$$U_{G,\alpha}(y)|_{y=-G^{-1}(1)} = G^{-1}(1) + \gamma_1(0) = G^{-1}(1) \geq 0,$$

with equality iff $G(0) = 1$. From these properties, we infer that

(i) $G(0) \in (0,1) \Rightarrow G \in \mathcal{C}_1$;
(ii) $G \in \mathcal{C}_2 \Rightarrow G(0) = 0$;
(iii) $G \in \mathcal{C}_3 \Rightarrow G(0) = 1$.

Furthermore, if $G(0) = 0$, then $\gamma_1 = 0$ and for such $G$'s: $U_{G,\alpha}(y) = -y \geq 0$ for all values $y \leq -G^{-1}(0) = 0$, implying that

(iv) $G(0) = 0 \Rightarrow G \in \mathcal{C}_2$.

Also, if $G(0) = 1$, then $-\gamma_1(1-\alpha) = G^{-1}(\alpha)$ and $G(-\gamma_1((1-\alpha)(1-G(-y)))) = \alpha(1-G(-y)) + G(-y) \geq G(-y)$; telling us that $U_{G,\alpha}(y) \leq 0$ for all $y \geq -G^{-1}(1)$, and implying that

(v) $G(0) = 1 \Rightarrow G \in \mathcal{C}_3$.

Finally the converse of (i) follows from (iv) and (v).

Although $y_0$ depends on $(\alpha, G)$, we will not stress this dependence unless necessary. Here are some useful facts concerning Definition 4's $y_0$. □

**Lemma 4.** (a) *For $G \in \mathcal{C}_1$, we have $U_{G,\alpha}(y) < 0$ iff $y > y_0$;*
(b) *Furthemore, we have*

(11) $$\gamma_2((1-\alpha)(1-G(-y_0))) = G^{-1}((1-\alpha) + \alpha G(-y_0)).$$

*Proof.* **(a)** We prove the result for $y > y_0$ only, with a proof for $y < y_0$ following along the same lines. We want to show that $U_{G,\alpha}(y) < 0$ for $y > y_0$, i.e.,

(12) $$-y < -\gamma_1((1-\alpha)(1-G(-y))).$$

Define

$$A = G(-y_0) - G(-y)$$
$$B_1 = G(-y_0) - G(-\gamma_1((1-\alpha)(1-G(-y)))) \text{ and}$$
$$B_2 = G(\gamma_2((1-\alpha)(1-G(-y)))) - G(\gamma_2((1-\alpha)(1-G(-y_0)))).$$

Observe that $A > 0$ since $-y < -y_0$. Since the quantities $\gamma_2(z), 1-G(-z)$, and $G(z)$ are all increasing in $z$, it follows as well that $B_2 \geq 0$. Now, with the definition of $\gamma_1$ and $\gamma_2$, and the identity $U_{G,\alpha}(y_0) = 0$, we have

$$B_1 + B_2 = (1-\alpha)(1-G(-y)) + G(-y_0) - G(\gamma_2((1-\alpha)(1-G(-y_0))))$$
$$= (1-\alpha)(1-G(-y))$$
$$\quad + G(-\gamma_1((1-\alpha)(1-G(-y_0)))) - G(\gamma_2((1-\alpha)(1-G(-y_0))))$$
$$= (1-\alpha)(1-G(-y)) - (1-\alpha)(1-G(-y_0))$$
$$= (1-\alpha)(G(-y_0) - G(-y)) = (1-\alpha)A < A \text{ (since } A > 0\text{).}$$

Finally, the inequality $B_1 < A$ is equivalent to (12) and establishes part (a) for $y > y_0$.

**(b)** Using the identity $U_{G,\alpha}(y_0) = 0$ and Definition 2, we have directly

$$-y_0 = -\gamma_1((1-\alpha)(1-G(-y_0)))$$
$$\Leftrightarrow \quad G(-y_0) = G(-\gamma_1((1-\alpha)(1-G(-y_0))))$$
$$\Leftrightarrow \quad G(-y_0) = G(\gamma_2((1-\alpha)(1-G(-y_0)))) - (1-\alpha)(1-G(-y_0))$$
$$\Leftrightarrow \quad 1 - \alpha + \alpha G(-y_0) = G(\gamma_2((1-\alpha)(1-G(-y_0))));$$

which is indeed equivalent to (11). □

**Lemma 5.** *Under Assumption 1,*

(a) *the $(1-\alpha) \times 100\%$ HPD credible interval $I_{\pi_0}(X)$ is of the form $[l(X), u(X)]$ with either:*

(13) $\quad$ (i) $l(x) = 0, u(x) = a_1(x) + a_2(x)G^{-1}(1 - \alpha + \alpha G(-\frac{a_1(x)}{a_2(x)}));$

*or*

(14)
$$\text{(ii)} \quad l(x) = a_1(x) - a_2(x)\gamma_1((1-\alpha)(1-G(-\frac{a_1(x)}{a_2(x)}))),$$
$$u(x) = a_1(x) + a_2(x)\gamma_2((1-\alpha)(1-G(-\frac{a_1(x)}{a_2(x)})));$$

(b) *Furthermore,* (i) *occurs iff* $G \in \mathcal{C}_2$ *or* $G \in \mathcal{C}_1$ *with* $\frac{a_1(x)}{a_2(x)} \leq y_0$; *(and equivalently* (ii) *occurs iff* $G \in \mathcal{C}_3$, *or* $G \in \mathcal{C}_1$ *with* $\frac{a_1(x)}{a_2(x)} \geq y_0$*).*

*Proof.* Part (b) follows from (14), the definition of the classes $\mathcal{C}_i; i = 1, 2, 3$; and Lemma 4. To establish part (a), proceed as in the proof of Theorem 1 by denoting $H_x$, and $H_x^{-1}$ as the cdf, and inverse cdf of the posterior distribution of $\tau(\theta)$ under $\pi_0$. Since $-T(X, \theta)$ is a pivot with cdf $G$, implying that its distribution for any given $\theta$ is free of $\theta$, we infer from (3) that, for $\theta \sim \pi$, $P_\pi(T(X, \theta) \geq y|x) = G(-y)$, or equivalently $P_\pi(\tau(\theta) \leq y|x) = G(\frac{y - a_1(x)}{a_2(x)})$. By definition of $\pi_0$, this gives us for $y \geq 0$,

$$H_x(y) = P_{\pi_0}(\tau(\theta) \leq y|x) = \frac{P_\pi(0 \leq \tau(\theta) \leq y|x)}{P_\pi(\tau(\theta) \geq 0|x)} = \frac{G(\frac{y-a_1(x)}{a_2(x)}) - G(-\frac{a_1(x)}{a_2(x)})}{1 - G(-\frac{a_1(x)}{a_2(x)})},$$

and

$$H_x^{-1}(\Delta) = a_1(x) + a_2(x)G^{-1}(\Delta + (1-\Delta)G(-\frac{a_1(x)}{a_2(x)})).$$

Now, observe that the posterior density $\tau(\theta)|x$ $(\propto G'(\frac{y-a_1(x)}{a_2(x)})I_{[0,\infty)}(y))$ is unimodal, with a maximum at $\max(0, a_1(x))$. Hence, we must have either: **(i)** $l(x) = 0, u(x) = H_x^{-1}(1-\alpha)$ yielding (13); or **(ii)**

$$H_x(u(x)) - H_x(l(x)) = 1 - \alpha, \text{ with } u(x) - l(x) \text{ minimal}$$
$$\iff G(\frac{u(x) - a_1(x)}{a_2(x)}) - G(\frac{l(x) - a_1(x)}{a_2(x)}) = (1-\alpha)(1 - G(-\frac{a_1(x)}{a_2(x)})),$$

with $u(x) - l(x)$ minimal,

yielding (14) by definition $\gamma_1$ and $\gamma_2$ (see Definition 2). □

**Theorem 2.** *Under Assumption 1, we have*

(a) *For $G \in \mathcal{C}_2$, $P_\theta(I_{\pi_0}(X) \ni \tau(\theta)) > 1 - \alpha$ for all $\theta$ such that $\tau(\theta) \geq 0$; and $P_\theta(I_{\pi_0}(X) \ni 0) = 1$ for all $\theta$ such that $\tau(\theta) = 0$;*

(b) *For $G \in \mathcal{C}_3$, $P_\theta(I_{\pi_0}(X) \ni \tau(\theta)) < 1 - \alpha$ for all $\theta$ such that $\tau(\theta) \geq 0$; and $P_\theta(I_{\pi_0}(X) \ni 0) = 0$ for all $\theta$ such that $\tau(\theta) = 0$;*

(c) *For $G \in \mathcal{C}_1$, $P_\theta(I_{\pi_0}(X) \ni \tau(\theta)) > (1-\alpha)(1 - G(-y_0))$ for all $\theta$ such that $\tau(\theta) \geq 0$; and $P_\theta(I_{\pi_0}(X) \ni 0) = 1 - G(-y_0)$ for all $\theta$ such that $\tau(\theta) = 0$ (with $y_0$ given in Definition 4);*

(d) *For unimodal $G'$, we have $\lim_{\tau(\theta) \to \infty} P_\theta(I_{\pi_0}(X) \ni \tau(\theta)) = 1 - \alpha$.*

*Proof.* **(a)** If $G \in \mathcal{C}_2$, $I_{\pi_0}(X)$ is given by (13) with probability one. This implies that $P_\theta(I_{\pi_0}(X) \ni 0) = 1$ for all $\theta$. As well, $u(x) \geq a_1(x) + a_2(x)G^{-1}(1-\alpha)$ implying that

$$\begin{aligned} P_\theta(I_{\pi_0}(X) \ni \tau(\theta)) &= P_\theta(\tau(\theta) \leq u(X)) \\ &\geq P_\theta(\tau(\theta) \leq a_1(X) + a_2(X)G^{-1}(1-\alpha)) \\ &= P_\theta(\frac{\tau(\theta) - a_1(X)}{a_2(X)} \leq G^{-1}(1-\alpha)) = G(G^{-1}(1-\alpha)) = 1 - \alpha. \end{aligned}$$

**(b)** If $G \in \mathcal{C}_3$, $I_{\pi_0}(X)$ is given by (14) with probability one. This implies that $P_\theta(I_{\pi_0}(X) \ni 0) = 0$ for all $\theta$ (in particular for those $\theta$ such that $\tau(\theta) = 0$). As well, since $l(x) = a_1(x) - a_2(x)\gamma_1((1-\alpha)(1 - G(-\frac{a_1(x)}{a_2(x)}))) \geq l(x) = a_1(x) - a_2(x)\gamma_1(1-\alpha)$, and similarly $u(x) \leq a_1(x) + a_2(x)\gamma_2(1-\alpha)$, we infer that

$$\begin{aligned} P_\theta(I_{\pi_0}(X) \ni \tau(\theta)) &\leq P_\theta(a_1(X) - a_2(x)\gamma_1(1-\alpha) \\ &\leq \tau(\theta) \leq a_1(X) + a_2(X)\gamma_2(1-\alpha)) \\ &= P_\theta(-\gamma(1-\alpha) \leq \frac{\tau(\theta) - a_1(X)}{a_2(X)} \leq \gamma_2(1-\alpha)) = 1 - \alpha. \end{aligned}$$

**(c)** First, given that coverage at $\tau(\theta) = 0$ occurs if and only if $l(X) = 0$, we have for $\theta$ such that $\tau(\theta) = 0$

$$P_\theta(I_{\pi_0}(X) \ni 0) = P_\theta(l(X) = 0) = P_\theta(\frac{a_1(X) - 0}{a_2(X)} \leq y_0) = 1 - G(-y_0).$$

For the more general lower bound, the idea here is the same as the one in Theorem 1, namely to work with a subset (with probability one) $I'(X)$ of $I_{\pi_0}(X)$ for which the coverage of $I'(X)$ is equal to $(1 - \alpha)(1 - G(-y_0))$. To achieve this, we first establish that

(15) $$u(x) \geq a_1(x) + a_2(x)\gamma_2((1-\alpha)(1 - G(-y_0))).$$

Indeed, if $\frac{a_1(x)}{a_2(x)} \leq y_0$, then $u(x) \geq a_1(x) + a_2(x)G^{-1}((1 - \alpha + \alpha G(-y_0)) = a_1(x) + a_2(x)\gamma_2((1-\alpha)(1-G(-y_0)))$, using (11). On the other hand, if $\frac{a_1(x)}{a_2(x)} \geq y_0$, then (15) follows directly as both $\gamma_2(z)$ and $1 - G(-z)$ increase with $z$. Similarly, if $\frac{a_1(x)}{a_2(x)} \geq y_0$, $l(x)$ is bounded above by $a_1(x) - a_2(x)\gamma_1((1-\alpha)(1-G(-y_0)))$. The above bounds on $l(x)$ and $u(x)$ imply that the coverage probability of $I_{\pi_0}(X)$ is bounded below by the coverage probability of $[\max(0, a_1(x) - a_2(x)\gamma_1((1-\alpha)(1 - G(-y_0))), a_1(x) + a_2(x)\gamma_2((1-\alpha)(1-G(-y_0)))]$; or equivalently by the coverage probability of

$$I'(X) = [a_1(x) - a_2(x)\gamma_1((1-\alpha)(1 - G(-y_0))), a_1(x) + a_2(x)\gamma_2((1-\alpha)(1-G(-y_0)))].$$

But finally, using the definition of $\gamma_1$ and $\gamma_2$, assumption (3), and the fact that $-T(X,\theta)|\theta$ has cdf $G$, we have $P_\theta(I_{\pi_0}(X) \ni \tau(\theta)) > P_\theta(I'(X) \ni \tau(\theta))$, with $P_\theta(I'(X) \ni \tau(\theta)) = P_\theta(a_1(X) - a_2(X)\gamma_1((1-\alpha)(1-G(-y_0))) \leq \tau(\theta) \leq a_1(X) + a_2(X)\gamma_2((1-\alpha)(1-G(-y_0)))) = P_\theta(-\gamma_1((1-\alpha)(1-G(-y_0))) \leq \frac{\tau(\theta)-a_1(X)}{a_2(X)} \leq \gamma_2((1-\alpha)(1-G(-y_0)))) = (1-\alpha)(1-G(-y_0))$.

**(d)** The result may be established along the lines of part (d) of Theorem 1. □

The next result is a specialization of Theorem 2 to cases where the density $G'$ is skewed to the right, or in other words the density of the pivot $\frac{a_1(X)-\tau(\theta)}{a_2(X)}$ is skewed to the left. Namely, the following result demonstrates that the lower bounds on the frequentist coverage probability for symmetric densities $G'$ also apply necessarily to right-skewed densities $G'$.

**Corollary 2.** *Under the conditions of Theorem 2 suppose further that*

(16) $$G(-\gamma_1(1-z)) \leq \frac{z}{2}; \text{ for all } z \in (0,1];$$

*(or equivalently $G(-\gamma_1(1-z)) \leq 1 - G(\gamma_2(1-z))$ by definition of $\gamma_1$ and $\gamma_2$). Then, under the assumptions of Theorem 2, we have*

(a) $P_\theta(I_{\pi_0}(X) \ni \tau(\theta)) > \frac{1-\alpha}{1+\alpha}$ *for all $\theta$ such that $\tau(\theta) \geq 0$;*
(b) $P_\theta(I_{\pi_0}(X) \ni 0) \geq \frac{1}{1+\alpha}$ *for all $\theta$ such that $\tau(\theta) = 0$.*

*Proof.* If $G \in \mathcal{C}_2$, the lower bounds hold of course by virtue of Theorem 2. If $G \in \mathcal{C}_1$, then we must have

$$G(-y_0) = G(-\gamma_1(1-\alpha)(1-G(-y_0))) \leq \frac{1-(1-\alpha)(1-G(-y_0))}{2},$$

with the inequality following from (16). The above now tells us that $1 - G(-y_0) \geq \frac{1}{1+\alpha}$, and parts (a) and (b) follow from Theorem 2. There remains to show that $G \in \mathcal{C}_3$ is incompatible with condition (16). But, if $G \in \mathcal{C}_3$ (i.e., $G(0) = 1$), then (16) cannot hold for $z = 1$ as $G(-\gamma_1(0)) = G(0) = 1 > \frac{1}{2}$. □

**Remark 2.** Corollary 2 includes the particular case of symmetry with equality in (16) and, therefore, can be viewed as a generalization of the results of Theorem 1. The lower bounds on coverage probability given in Theorem 2 and Corollary 2 correspond to the ones given by Zhang and Woodroofe [8] for a lower bounded scale parameter of a Fisher distribution, and arising in the estimation of the ratio of variance components in a one-way balanced model analysis of variance with random effects.

**Example 2** (Lower bounded Gamma scale parameter). As a followup to Example (e) of Section 3, consider a Gamma$(r,\theta)$, $\theta \geq a > 0$, where $X|\theta \sim \frac{1}{\theta}f_1(\frac{x}{\theta})$ with $f_1(y) = \frac{y^{r-1}e^{-y}}{\Gamma(r)}1_{(0,\infty)}(y)$. Considering the cdf $G$ of $-(\log(\frac{X}{\theta}) - m)$, where $m$ is chosen in such a way that $G'$ has a mode at 0, which is required in Lemma 5 and Theorem 2, we obtain that $-(\log(\frac{X}{\theta}) - \log r)$ has pdf

$$G'(y) = \frac{r^r}{\Gamma(r)}e^{-r(y+e^{-y})},$$

and cdf $G(y) = P(\text{Gamma}(r,1) \geq re^{-y})$. For instance, with the exponential case ($r = 1$), we have $G(y) = e^{-e^{-y}}$. Definition 2's $\gamma_1(\Delta)$ and $\gamma_2(\Delta)$ satisfy the equation $-\gamma_1(\Delta) + e^{\gamma_1(\Delta)} = \gamma_2(\Delta) + e^{-\gamma_2(\Delta)}$ with $G(\gamma_2(\Delta)) - G(-\gamma_1(\Delta)) = \Delta$, but are not available explicitly. Hence, as will be the case in general, neither the lower and upper bounds $l(X)$ and $u(X)$ of Lemma 5, nor Definition 4's $y_0$, are available explicitly. However, Theorem 2 (parts c and d) do apply. For instance, with $r = 3, \alpha = 0.05$, a numerical evaluation yields $y_0 \approx 0.912968$ and $1 - G(-y_0) \approx 0,979353$, which gives the exact coverage at the boundary $\theta = a$, and which tells us that Theorem 2's lower bound on coverage $(1-\alpha)(1-G(-y_0))$ is approximatively equal to $(0.95)(0.979353) = 0.930386$. We were unable to establish but believe that (16) holds for the cdf's of this example, which would permit the application of Corollary 2, but observe that the lower bound of 0.930386 actually exceeds Corollary 2's lower bound of $\frac{0.95}{1.05} = 0.904762$.

We conclude by pointing out that the results of this paper do leave open several questions concerning further coverage probability properties of the Bayesian confidence interval $I_{\pi_0}(X)$. Namely, as seen in the above example, it would be desirable for the quantity $1 - G(-y_0)$ of Theorem 2 to be made more explicit. Further numerical evaluations of $1 - G(-y_0)$, which also suggest quite high lower bounds on coverage, are given by Zhang and Woodroofe [8] in their particular case of a lower bounded Fisher distribution scale parameter.

## Acknowledgments

The authors thank Amir Payandeh and Keven Bosa for a diligent reading and useful comments. The authors are thankful to Claude Beauchamp who assisted with numerical evaluations. For the work of Marchand, the support of NSERC of Canada is gratefully acknowledged, while for the work of Strawderman, the support of NSA Grant 03G-1R is gratefully acknowledged. Finally, we are grateful to an anonymous referee and Associate Editor for constructive comments that led to an improved manuscript.